\documentclass[12pt]{amsart}
\newcommand{\RP}{{\mathbb {RP}}}
\newcommand{\CP}{{\mathbb {CP}}}
\newcommand{\Z}{{\mathbb Z}}
\newcommand{\R}{{\mathbb R}}

\newtheorem{example}{Example}
\newtheorem{proposition}{Proposition}

\newtheorem{definition}{Definition}
\newtheorem{lemma}{Lemma}

\begin{document}
\author{V.A.~Vassiliev}
\address{Steklov Mathematical Institute of Russian Academy of Sciences;
\newline \hspace*{3mm}
National Research University Higher School of Economics} \email{vva@mi.ras.ru}


\dedicatory{To the memory of S.~Duzhin}

\title{Linking numbers in non-orientable 3-manifolds}

\begin{abstract}
The construction of integer linking numbers of closed curves in a
three-dimensional manifold usually appeals to  the orientation of this
manifold. We discuss how to avoid it constructing similar homotopy invariants
of links in non-orientable manifolds.
\end{abstract}

\maketitle

\section{Introduction}
Let $M^3$ be a connected three-dimensional manifold, and ${\mathcal C}$ the
disjoint union of two circles, ${\mathcal C}=S^1_1 \sqcup S^1_2$. A smooth map
$f:{\mathcal C} \to M^3$ is a {\em link} if it is a smooth embedding. Two links
are {\em homotopy equivalent} (see \cite{Milnor}, \cite{Mel}) if there is a
homotopy $F: {\mathcal C} \times [0,1] \to M^3$ between them, such that for any
$\tau \in [0,1]$ the sets $F(S^1_1 \times \tau)$ and $F(S^1_2 \times \tau)$
have no common points. If $M^3$ is orientable and both parts $f(S^1_1),
f(S^1_2)$ of the link $f$ are contractible in $M^3$, then their {\em linking
number} is an integer-valued invariant separating some classes of homotopy
equivalence of such links. It is equal to 0 if these two parts lie in two
non-intersecting embedded balls in $M^3$; for two links which differ only in a
small domain and look there as shown in two sides of Fig. \ref{fig1} the
linking number of the left-hand link is greater by 1 than that of the
right-hand one if the orientation of the ambient manifold is as shown in the
center of the picture, and is smaller by 1 otherwise. This invariant can be
defined equivalently as the intersection number of the path connecting our
given link with the distinguished one in the space of pairs of curves in $M^3$,
and the {\em discriminant} set in this space, which consists of maps $f$ with
$f(S^1_1) \cap f(S^1_2) \neq \emptyset$; the coorientation of this set is
derived from that of $M^3$ by the rule of Fig.~\ref{fig1}.

\unitlength=0.4mm \linethickness{0.4pt}
\begin{figure}
\begin{picture}(120,33)
\put(0,3){\vector(1,1){30}} \put(30,3){\line(-1,1){14}}
\put(14,19){\vector(-1,1){14}} \put(12,0){\large +} \put(90,3){\line(1,1){14}}
\put(106,19){\vector(1,1){14}} \put(120,3){\vector(-1,1){30}}
\put(102,0){\large $-$} \put(60,15){\vector(1,0){13}}
\put(60,15){\vector(0,1){13}} \put(60,15){\vector(-1,-1){10}} \put(68,9){\tiny
$1$} \put(55,23){\tiny $2$} \put(50.5,9.5){\tiny $3$}
\end{picture}
\caption{Coorientation of the discriminant variety} \label{fig1}
\end{figure}
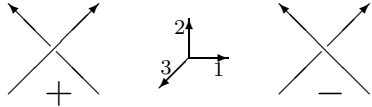

\unitlength=0.8mm \linethickness{0.4pt}

More generally, one can try to separate links within one connected component of
the space $C^\infty({\mathcal C},M^3)$ by a similar invariant in the case of
links with non-contractible parts, see especially the work \cite{CR} in which
the multidimensional case is also studied. Again, the linking number is set to
be equal to 0 at some distinguished link from this connected component, and is
extended by the rule of Fig. \ref{fig1} to all other links in it. However, this
definition is self-consistent only under some additional conditions; on the
other hand, the obtained invariant takes values in a group which can be greater
than just $\Z$.

Sometimes there exists a canonical way to define such a distinguished link. For
instance, it is so if we consider the links, only one part $f(S^1_i)$ of which
is contractible in $M^3$: this part of the distinguished link should be
contained in a ball separated from the other part. Also, if $M^3 = M^2 \times
\R^1$ then we can choose a link $f$ in any component of $C^\infty({\mathcal
C},M^3)$ in such a way that $f(S^1_1) \subset M^2 \times \R^1_-$ and $f(S^1_2)
\subset M^2 \times \R^1_+$, cf. \cite{GVsplit}.

In any way, the orientation of the ambient manifold plays the crucial role in
the construction of integer linking numbers; its direct extension to the
non-orientable case can provide the mod 2 indices only. We show below how to
improve it to define integer homotopy invariants for links in non-orientable
manifolds. All these invariants will be defined in the terms of surgeries
similar to the one shown in Fig.~\ref{fig1}. So, they can be considered as
generalizations of the usual linking number, and will be called the
(generalized) linking numbers throughout this article.

Some related questions and examples for the similar problem of isotopy
classification of knots were discussed in \cite{Vas97}, \cite{FiBook} and
\cite{GMV}. Almost all our examples concern the links in thickened surfaces. A
natural problem is to realize these invariants for links in such manifolds in
the terms of arrow diagrams in the spirit of \cite{FiArt}, \cite{PV},
\cite{FiBook}. However, the general constructions described below can be
applied to arbitrary 3-manifolds as well.

\section{First degree linking numbers}
\label{first}

\subsection{Finite-dimensional approximations}
\label{findim} We will work with the space $C^\infty({\mathcal C},M^3)$ as with
an open domain in a real affine space of a very large but finite dimension. To
justify this assumption, let us fix a smooth embedding $I: M^3 \to \R^n$ ($n$
large enough), a tubular neighborhood $U \subset \R^n$ of the submanifold
$I(M^3)$, and a smooth map $\psi: {\mathcal C} \to U$. Consider the space of
all maps ${\mathcal C} \to U$ of the form $\psi+ p_d$ where $p_d$ is defined by
a collection of $2n$ Fourier polynomials of a sufficiently large degree $d$.
Such maps composed with the canonical projection $U \to I(M^3) \simeq M^3$
provide a finite-dimensional approximation of the space $C^\infty({\mathcal
C},M^3)$; let us denote it by ${\mathcal F}_d$. Any link can be approximated by
such a map so that the approximating link is in the same homotopy class.
Moreover, if two links from ${\mathcal F}_d$ are in one and the same homotopy
equivalence class in $C^\infty({\mathcal C},M^3)$, then such a homotopy can be
approximated by a homotopy connecting these two links within some space
${\mathcal F}_{d'}$, maybe with $d'>d$. All topological properties of sets of
links and singular links from ${\mathcal F}_d$, which we will use, stabilize
when $d$ grows to infinity.

\subsection{The discriminant and its components}

The simplest invariant of a two-component link in $M^3$ is the corresponding
connected component of the space $C^\infty({\mathcal C},M^3) \simeq
(C^\infty(S^1,M^3))^2$. The set of such components is in the obvious one-to-one
correspondence with the set $(\pi_1(M^3)/\mbox{conj})^2$ where conj is the
equivalence relation identifying any element $\alpha \in \pi_1(M^3)$ with all
elements of type $c^{-1}\alpha c$. The homotopy classes within one and the same
component are separated by the discriminant set $\Sigma$, i.e. the set of links
whose two parts have common points. If $d$ is not too small, then this set is a
hypersurface in ${\mathcal F}_d$. It is convenient to assume that $M^3$, its
embedding $I: M^3 \hookrightarrow \R^n$, projection $U \to M^3$ and map $\psi:
{\mathcal C} \to U$ are analytic, then this hypersurface $\Sigma$ is
semianalytic. If $M^3$ is oriented, then $\Sigma$ is (co)oriented in its
non-singular points by the rule of Fig.~\ref{fig1}, and defines an integer
cycle with closed supports in any ${\mathcal F}_d$ with sufficiently large $d$.
The standard linking number of a link (if it is well-defined) is the
intersection number of this cycle and any path connecting our link with the
distinguished one; this number does not change when $d$ grows. Any irreducible
component of the discriminant in ${\mathcal F}_d$ belongs to the intersection
of ${\mathcal F}_d$ with some irreducible component of the discriminant in
${\mathcal F}_{d'}$ for any $d'>d$. By the set of irreducible components of the
discriminant in $C^\infty({\mathcal C}, M^3)$ we mean the direct limit of the
arising directed set. In what follows we will skip similar limit constructions,
speaking directly on subsets etc. in $C^\infty({\mathcal C}, M^3)$.

\begin{proposition} \label{prop1}

A. The irreducible components of the discriminant set in $C^\infty({\mathcal
C},M^3)$ are in a natural one-to-one correspondence with the elements of
$(\pi_1(M^3))^2/\mbox{conj}$, i.e. with the ordered pairs $(\alpha,\beta)$ of
elements of $\pi_1(M^3)$ considered up to the simultaneous conjugation:
$(\alpha,\beta) \sim (\alpha',\beta')$ if and only if there is $c\in
\pi_1(M^3)$ such that $\alpha'=c^{-1}\alpha c$, $\beta'=c^{-1}\beta c$.

B. An irreducible component of $\Sigma$ represented by such a pair
$(\alpha,\beta)$ is {\em not coorientable} in $C^\infty({\mathcal C},M^3)$ if
and only if there exists $c \in \pi_1(M^3)$ reversing the orientation of $M^3$
and commuting with both $\alpha$ and $\beta$. $($In particular, this condition
is never satisfied if $M^3$ is orientable$)$.
\end{proposition}

{\it Proof.} A. Define the {\em tautological resolution} $\sigma_1$ of the
discriminant $\Sigma \subset C^\infty({\mathcal C},M^3)$ as the submanifold of
codimension 3 in the space $S^1_1 \times S^1_2 \times C^\infty({\mathcal
C},M^3),$ consisting of triples $(x,y,f)$ such that $f(x)=f(y)$. Obviously, it
is a smooth submanifold in this space. The forgetful map $\sigma_1 \to \Sigma$
sending any triple $(x,y,f)$ to $f$ is surjective and is a local diffeomorphism
close to the {\em generic} discriminant maps $f$ having only one transverse
intersection point of $f(S^1_1)$ and $f(S^1_2)$. Irreducible components of
$\Sigma$ are the images of connected components of $\sigma_1$. Also, we have a
locally trivial fiber bundle $p:\sigma_1 \to M^3$ sending any triple $(x,y,f)$
to the point $f(x)\equiv f(y)$. The connected components of its fiber over a
point $m \in M^3$ are in the obvious one-to one correspondence with the
elements of $(\pi_1(M^3,m))^2$. We get the covering over $M^3$ with the fiber
$(\pi_1(M^3,m))^2$; the monodromy of this covering over the loop $c$ sends any
such component $(\alpha, \beta)$ to $(c^{-1}\alpha c, c^{-1}\beta c)$.

B. We can choose arbitrarily an orientation of $M^3$ inside a
simply-con\-nec\-ted neighborhood $O(m)$ of the marked point $m$. The
coorientation (and hence also orientation) of the piece of any irreducible
component of $\Sigma$ consisting of maps $f$ having intersection points in this
neighborhood only is well defined by the rule of Fig.~\ref{fig1}: the link
marked by $+$ (respectively, $-$) lies on the positive (respectively, negative)
side of the discriminant. This orientation induces an orientation of the
manifold $p^{-1}(O(m)) \subset \sigma_1$. This choice of orientations fails
when we go over a loop $c \in \pi_1(M^3,m)$ reversing the orientation of $M^3$.
However, the monodromy over such a loop $c$ leads us to the same component
$(\alpha, \beta)$ (and hence reverses its orientation) only if $c$ commutes
with both $\alpha$ and $\beta$. \hfill $\Box$

\subsection{Homological condition}
\label{homc} If some irreducible component of $\Sigma$ is cooriented, then the
integer intersection index of this component with any generic path in
$C^\infty({\mathcal C},M^3)$ connecting some two non-discriminant links is
well-defined. However, sometimes this index is not an invariant of the
endpoints of the path: it can happen that such intersection indices defined by
different paths with the same endpoints are different. This situation is
equivalent to the condition that the Poincar\'e intersection pairing with our
component of $\Sigma$ defines a non-trivial element of $H^1(C^\infty({\mathcal
C},M^3), \Z)$. Therefore we need to check whether this happens or not; of
course, it is enough to calculate the intersection numbers of our component of
$\Sigma$ with some collection of generating elements of $H_1(C^\infty({\mathcal
C},M^3), \Z)$ only.

\subsection{Examples} By Proposition \ref{prop1}(B), we have no chance to
gain a generalized linking number defined as above if $M^3$ is non-orientable
and both parts of $f({\mathcal C})$ are contractible. Also, we have no chance
if $M^3$ is equal to $\RP^2 \times M^1$, since in this case the group
$\pi_1(M^3)$ is commutative.

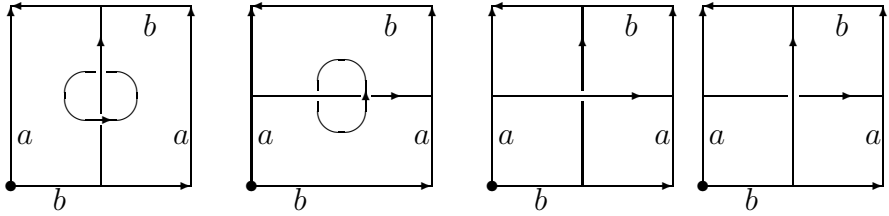
\begin{figure}
\begin{picture}(145,35)
\put(0,5){\circle*{1.7}} \put(0,5){\vector(1,0){30}}
\put(30,5){\vector(0,1){30}} \put(7,1){$b$} \put(27,12){$a$}
\put(30,35){\vector(-1,0){30}} \put(0,5){\vector(0,1){30}} \put(22,30){$b$}
\put(1,12){$a$} \put(15,5){\line(0,1){10}} \put(15,17){\line(0,1){18}}
\put(14,20){\oval(10,8)[l]} \put(16,20){\oval(10,8)[r]}
\put(14,16){\line(1,0){2}} \put(15,27){\vector(0,1){3}}
\put(13,16){\vector(1,0){4}}

\put(40,5){\circle*{1.7}} \put(40,5){\vector(1,0){30}}
\put(70,5){\vector(0,1){30}} \put(47,1){$b$} \put(67,12){$a$}
\put(70,35){\vector(-1,0){30}} \put(40,5){\vector(0,1){30}} \put(62,30){$b$}
\put(41,12){$a$} \put(40,20){\line(1,0){18}} \put(60,20){\line(1,0){10}}
\put(55,21){\oval(8,10)[t]} \put(55,19){\oval(8,10)[b]}
\put(59,19){\line(0,1){2}} \put(59,18){\vector(0,1){3}}
\put(61,20){\vector(1,0){4}}

\put(80,5){\circle*{1.7}} \put(80,5){\vector(1,0){30}}
\put(110,5){\vector(0,1){30}} \put(87,1){$b$} \put(107,12){$a$}
\put(110,35){\vector(-1,0){30}} \put(80,5){\vector(0,1){30}} \put(102,30){$b$}
\put(81,12){$a$} \put(80,20){\line(1,0){30}} \put(95,5){\line(0,1){14}}
\put(95,21){\line(0,1){14}} \put(95,27){\vector(0,1){3}}
\put(101,20){\vector(1,0){4}}

\put(115,5){\circle*{1.7}} \put(115,5){\vector(1,0){30}}
\put(145,5){\vector(0,1){30}} \put(122,1){$b$} \put(142,12){$a$}
\put(145,35){\vector(-1,0){30}} \put(115,5){\vector(0,1){30}} \put(137,30){$b$}
\put(116,12){$a$} \put(115,20){\line(1,0){14}} \put(131,20){\line(1,0){14}}
\put(130,5){\line(0,1){30}} \put(130,27){\vector(0,1){3}}
\put(136,20){\vector(1,0){4}}

\end{picture}
\caption{Examples of links in the thickened Klein bottle} \label{fig2}
\end{figure}

Consider the next complicated case $M^3=K^2 \times \R^1$, where $K^2$ is the
Klein bottle, see Fig.~\ref{fig2}. For the marked point $m \in K^2$ we choose
the one represented by all corners of the square shown in this picture. For the
fundamental domain of $K^2$ we take this square less its upper and right-hand
sides, in particular the canonical representative of $m$ is the lower bottom
corner. Any fundamental group $\pi_1(K^2,m')$, $m' \neq m$, will be identified
with $\pi_1(K^2,m)$ by means of the segment inside this fundamental domain
connecting these points $m$ and $m'$. The group $\pi_1(M^3,(m \times 0))\equiv
\pi_1(K^2,m)$ is generated by two elements $a, b$ with one relation $a=bab$.
Any element of this group can be reduced to the normal form $a^rb^s$, $r$ and
$s$ integers. Such an element reverses the orientation of $M^3$ if and only if
$r$ is odd.

In our first example (Fig.~\ref{fig2} left) one part of the link is
contractible, and the other one is homotopic to the loop $a$. A path in
$C^\infty({\mathcal C}, M^3)$ connecting this link with the distinguished one
crosses the discriminant at unique point of its component whose code $(\alpha,
\beta)$ in the sense of Proposition \ref{prop1} is equal to $(1, a)$. The
element $\{a\} \in \pi_1(K^2)$ commutes with both $1$ and $a$ and reverses the
orientation of $M^3$, hence this component is not coorientable and does not
define an integer linking number.
\medskip

In the next example (second from the left in Fig.~\ref{fig2}) $f(S^1_1)$ is
again a contractible knot, and $f(S^1_2)$ is homotopic to $b$ and to $b^{-1}$.
The unique component of $\Sigma$ separating our link from the distinguished one
with the same projection to $K^2$ has the code $(\alpha, \beta) = (1, b)$. It
is easy to calculate that the element $b \in \pi_1(M^3)$ commutes only with
elements $a^rb^s$ with even $r$, and these are exactly the elements that
preserve the orientation of $M^3$. Therefore the corresponding component of
$\sigma_1$ is orientable. Moreover, the homological condition of \S \ref{homc}
is satisfied by any of the following two statements.

\begin{proposition}
\label{prop55} If one of two parts of the link $f({\mathcal C})$ is
contractible and $\pi_2(M^3)=0$, then the intersection index with any
irreducible component of $\Sigma$ $($and, moreover, with any cycle in
$C^\infty({\mathcal C},M^3)$ contained in $\Sigma)$ defines the trivial element
of the one-dimensional cohomology group of the connected component of
$C^\infty({\mathcal C},M^3)$ containing the link $f$.
\end{proposition}

{\it Proof}. In this case any 1-dimensional homology class in this component of
$C^\infty({\mathcal C},M^3)$ can be realized by a one-dimensional family of
links (parameterized by the points of a circle) not intersecting the
discriminant at all: the contractible component of this link can be kept within
a small ball not meeting the second component. \hfill $\Box$ \medskip

\begin{proposition}
\label{prop66} If $M^3 = M^2 \times \R^1$ then the group
$H_1(C^\infty({\mathcal C}, M^3))$ is generated by the classes of loops lying
in the space of non-discriminant links.
\end{proposition}

Indeed, in this case we can realise any 1-dimensional homology class of
$C^\infty({\mathcal C}, M^3)$ by a one-parameter family of maps $f_\tau$ such
that $f_\tau (S^1_1)$ always belongs to $M^2 \times \R^1_-$, and $f_\tau
(S^1_2)$ to $M^3 \times \R^1_+$. \hfill $\Box$
\medskip

So, the component of $\Sigma$ containing the element $(1,b)$ defines well a
generalized linking number for links in $K^2 \times \R^1$, in particular it
separates the link shown in the second from the left picture of Fig.~\ref{fig2}
from the trivial one.
\medskip

Further, consider two links shown in the right-hand part of Fig. \ref{fig2}.
The code $(\alpha,\beta)$ for the discriminant curve $f({\mathcal C})$
separating them is equal to $(a,b)$ (which is equivalent to $(a,b^{-1})$ by
means of the simultaneous conjugation by $a$). The corresponding component of
$\Sigma$ is coorientable by the same reason as in the previous example, so that
this component is an integer cycle. Its dual cohomology class is equal to 0 by
Proposition \ref{prop66}, therefore we obtain an integer linking number taking
different values on our two links.
\medskip

In the last example $M^3 =S^2 \times S^1$. Consider the component of
$C^\infty({\mathcal C},M^3)$ containing the maps $f$ such that $f(S^1_1)$ is
contractible, and $f(S^1_2)$ runs once the circle $s \times S^1$, where $s$ is
a fixed point of $S^2$. The discriminant in this component consists of only one
irreducible component since $\pi_1(S^1 \times S^2)$ is commutative. The latter
component is orientable since $M^3$ is. However, it does not define a linking
number. Indeed, let the part $f(S^1_2)$ of our link stay unmoved, and
$f(S^1_1)$ run along an 1-cycle generating the group $\pi_1(\Omega(S^2)) \equiv
\pi_2(S^2) \simeq \Z$. The obtained 1-cycle in $C^\infty({\mathcal C},M^3)$ has
non-zero intersection index with the basic cycle of the discriminant.

\section{Second degree linking numbers of trivial knots}

As was mentioned above, there is no chance to define integer linking numbers
(i.e. the homotopy invariants Alexander dual to appropriate irreducible
components of the discriminant) separating links consisting of two contractible
knots in a non-orientable manifold: all such components are non-orientable.
However, sometimes it is possible to define similar invariants of higher
orders, separating such links, if we take in account the intersections and
self-intersections of these components. (In a similar way, the Boy surface is
not orientable, still it separates some points of its complement in $\R^3$
because it has self-intersections). The construction of these invariants is
almost the same as in \cite{knsp}, \cite{Vas97}: we need only to take
additional care on the (co)orientability of strata of the discriminant (holding
automatically in the case of oriented 3-manifolds). This construction is based
on the simplicial resolution of the discriminant constructed in the following
way (generalizing the constructions of \cite{knsp} and \cite{Gorinov}).

\subsection{The simplicial resolution and its filtration} \label{simpl}

For any finite unordered collection $J = \langle j_1, \dots, j_k\rangle$ of
natural numbers, each of which is greater than $1$, consider the configuration
space $Y_J$ of all finite unordered collections $\langle A_1, \dots, A_k
\rangle$ of disjoints subsets $A_i \subset {\mathcal C}$ of cardinalities
$j_i$, such that each of these subsets $A_i$ contains points of both $S^1_1$
and $S^1_2$. The union of such configuration spaces over all $J$ is supplied
with the following version of the Hausdorff metric: the distance between
configurations $\langle A_1, \dots, A_k \rangle$ and $\langle B_1, \dots, B_l
\rangle$ is equal to
\begin{equation}
\max_{i=1}^k \min_{i'=1}^l \mbox{dist}(A_i, B_{i'}) + \max_{i'=1}^l
\min_{i=1}^k \mbox{dist}(A_i, B_{i'}), \label{haus2}
\end{equation} where
$\mbox{dist}(A_i,B_{i'})$ is the usual Hausdorff distance
\begin{equation}
\max_{x\in A_i} \min_{y\in B_{i'}} \rho(x, y) + \max_{y\in B_{i'}} \min_{x\in
A_i} \rho(x, y),
\end{equation}
and $\rho(x,y)$ is equal to the angular distance in $S^1$ (taking values in
$[0,\pi]$) if $x$ and $y$ are in the same part $S^1_1$ or $S^1_2$ of ${\mathcal
C}$ and is equal to an arbitrary fixed number greater than $\pi$ if $x$ and $y$
are in different parts.

\begin{definition} \rm
The {\em complexity} of a collection $J = \langle j_1, \dots, j_k\rangle$ is
the number $\sum_{i=1}^k (j_i -1)$. For any natural $p$ denote by $Y(p)$ the
union of configuration spaces $Y_J$ over all collections $J$ of complexity
$\leq p$. \end{definition}

In particular, for any configuration $\langle A_1, \dots, A_k \rangle \in Y_J$
the codimension in $C^\infty({\mathcal  C}, \R)$ of the space of functions
taking equal values on all points any set $A_i$ is equal to the complexity of
$J$.

\begin{lemma}
A $($cf. \cite{Gorinov}$)$. Any space $Y(p)$ supplied with the metric
$($\ref{haus2}$)$ is compact.

B. The union of all spaces $Y(p)$ can be supplied with the structure of a
$CW$-complex such that any $Y(p)$ is its finite subcomplex; the $CW$-topology
of any of these subcomplexes is equivalent to that defined by the metric
$($\ref{haus2}$)$ $($although the $CW$-topology of entire this union does not
coincide with the one given by this metric$)$.
\end{lemma}

{\em Proof} is elementary. \hfill $\Box$ \medskip

Denote by $Y$ this union of spaces $Y(p)$ with this structure of a
$CW$-complex. Define a partial order on it: a configuration $\langle A_1, \dots
, A_k \rangle$ is subordinate to $\langle B_1, \dots, B_l \rangle$ if for any
$i=1, \dots , k$ there is $i' \in \{1, \dots, l\}$ such that $A_i \subset
B_{i'}$.

The $k$th {\em self-join} $Y^{*k}$ of $Y$ is defined as follows. For any
fixed $p$ consider a {\em generic} embedding $\varphi: Y(p) \hookrightarrow
\R^W$ into an Euclidean space of a very large dimension, and take the union
$(Y(p))^{*k}$ of all $(k-1)$-dimensional simplices in $\R^W$ all vertices of
which belong to $\varphi(Y(p))$. The genericity condition of the embedding
$\varphi$ means that these simplices do not have unexpected intersections: the
intersection set of any two simplices is some their common face. Such
embeddings are present in the space of all continuous maps $Y(p) \to \R^W$ if
$W$ is sufficiently large with respect to $p$ and $k$. Moreover, the spaces
$(Y(p))^{*k}$ defined in this way by different generic embeddings into (maybe)
different spaces $\R^W$ are canonically homeomorphic to one another and
canonically include all spaces $(Y(p'))^{*k'}$ with $p' \leq p$ and $k' \leq
k$. This allows us to define the space $Y^{*\infty}$ as the union of all these
spaces $(Y(p))^{*k}$ embedded in one another, with the direct limit topology.

\begin{definition} \rm
A $(k-1)$-dimensional simplex in $(Y(p))^{*k}$ with vertices in $\varphi(Y(p))$
is called {\em coherent} if the configurations corresponding to its vertices
are all incident to one another (and hence constitute a monotone sequence) in
the sense of the above defined partial order on the set $Y$. The space $\Xi(k)$
is defined as the union of all coherent $r$-dimensional simplices in
$(Y(k))^{*k}$, $r \leq k-1$; in particular $\Xi(k)$ is canonically embedded
into $\Xi(k')$, $k' > k$. The subspace $\Xi \subset Y^{*\infty}$ is defined as
the union of all these spaces $\Xi(k)$.
\end{definition}

The simplicial resolution $\sigma$ of the discriminant variety $\Sigma \subset
C^\infty({\mathcal C}, M^3)$ will be constructed as a subspace in
$C^\infty({\mathcal C}, M^3) \times \Xi$.

\begin{proposition}
For any sufficiently large natural numbers $n$ and $d$ there is an open dense
subset in the space $C^\infty({\mathcal C}, U)$ such that if the map $\psi$
participating in the definition of the space ${\mathcal F}_d$ $($see \S
\ref{findim}$)$ belongs to this subset, then for any $f \in {\mathcal F}_d$ the
set $f(S^1_1) \cap f(S^1_2)$ is finite, the complete preimage
\begin{equation}
\label{conff} \langle f^{-1}(m_1), \dots, f^{-1}(m_k) \rangle
\end{equation} of this set also
is finite, and the complexity of this preimage does not exceed the dimension of
${\mathcal F}_d$.
\end{proposition}

{\it Proof} follows almost immediately from the Thom's multijet transversality
theorem, see e.g. \cite{GG}. \hfill $\Box$ \medskip

We will assume that the conditions of this proposition are satisfied, in
particular for any $f \in \Sigma$ this preimage (\ref{conff}) is a
configuration defining a point of some space $Y(p)$, $p \leq \dim {\mathcal
F}_d$. Define the subcomplex $\Xi(f) \subset \Xi_p$ as the union of coherent
simplices all whose vertices correspond to the configurations subordinate to
(\ref{conff}). It is a finite simplicial complex, all simplices of which have a
common vertex; in particular it is contractible. The resolution $\sigma$ of
$\Sigma$ is defined as the union of complexes $f \times \Xi(f) \subset \Sigma
\times \Xi$ over all $f \in \Sigma$. Its obvious projection to $\Sigma \subset
{\mathcal F}_d$ is a surjective proper map with compact contractible preimages
of all points, in particular it induces an isomorphism of homology groups of
locally finite singular chains in $\sigma$ and $\Sigma$.

The space $\sigma$ has a natural increasing filtration: its term $\sigma_p$ is
the preimage of the subspace $\Xi(p) \subset \Xi$ under the obvious projection
of $\sigma \subset \Sigma \times \Xi$ to $\Xi$.

\begin{example} \rm The first term $\sigma_1$ of this filtration is exactly the
tautological resolution of $\sigma$ considered in \S \ref{first}. Indeed, the
unique collection $\langle j_1, \dots, j_k \rangle $ of complexity 1 with \ \
$\min j_i \geq 2$ \ \ is the one-member collection $\langle 2\rangle,$ related
with configurations of two points $\langle x \in S^1_1, y \in S^1_2\rangle $
and with the simplest chord diagram
\begin{picture}(12,5) \put(2,2){\circle{4}} \put(10,2){\circle{4}}
\put(4,2){\line(1,0){4}} \end{picture} . The collections $\langle j_1, \dots
j_k \rangle $ of complexity 2 generating the rest of the term $\sigma_2$ are
equal to $\langle 2,2\rangle $ and $\langle 3\rangle$; they are related with
chord diagrams \unitlength=1mm \linethickness{0.4pt}
\begin{picture}(12,5)
\put(2,2){\oval(4,4)} \put(10,2){\oval(4,4)} \put(4,1){\line(1,0){4}}
\put(4,3){\line(1,0){4}} \end{picture} and
\begin{picture}(12,5)
\put(2,2){\oval(4,4)} \put(10,2){\oval(4,4)} \put(4,1){\line(4,1){4}}
\put(4,3){\line(4,-1){4}} \end{picture} or
\begin{picture}(12,5)
\put(2,2){\oval(4,4)} \put(10,2){\oval(4,4)} \put(4,2){\line(4,1){4}}
\put(4,2){\line(4,-1){4}} \end{picture}.
\end{example}

\medskip

This filtration defines a spectral sequence, calculating the homology group of
locally finite chains in $\sigma$ (and hence also in $\Sigma$). Any its
homology class of maximal dimension (equal to the dimension of $\Sigma$)
defines the Poincar\'e dual element in the group $H^1(C^\infty({\mathcal
C},M^3),\Z)$. If this element is equal to zero, then our homology class defines
a homotopy invariant of links: a generalized linking number.

The formal description of the algorithm calculating these invariants (and also
of the stabilization of our spectral sequences for different approximating
spaces ${\mathcal F}_d$ with growing $d$) follows the one given in \cite{knsp}
and in many subsequent works. In our present case this algorithm has important
simplifications, because no analogs of the 1-term relation occur, and the
4-term relation is replaced just by the 2-term relation \unitlength=0.9mm
\linethickness{0.4pt}
\begin{picture}(15,5) \put(0,-1){\line(0,1){6}}
\put(4,-1){\line(0,1){6}} \put(0,0){\line(1,0){4}} \put(0,4){\line(1,0){4}}
\put(6,1.5){\tiny $=$} \put(11,-1){\line(0,1){6}} \put(15,-1){\line(0,1){6}}
\put(11,0){\line(1,1){4}} \put(11,4){\line(1,-1){4}} \end{picture} .
 On the other hand, one needs to take
into account the topology of the ambient manifold, including its
non-orientability and fundamental group. The next example demonstrates main
difficulties of this kind.

\unitlength=0.8mm \linethickness{0.4pt}
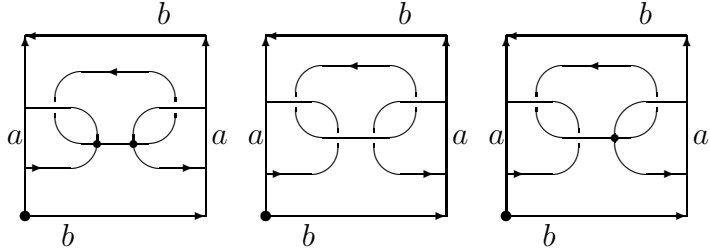
\begin{figure}
\begin{picture}(110,40)
\put(0,5){\circle*{1.7}} \put(0,13){\vector(1,0){3}}
\put(25,13){\vector(1,0){3}} \put(16.5,29){\vector(-1,0){3}}
\put(0,5){\vector(1,0){30}} \put(30,5){\vector(0,1){30}} \put(6,0){$b$}
\put(31,17){$a$} \put(30,35){\vector(-1,0){30}} \put(0,5){\vector(0,1){30}}
\put(22,37){$b$} \put(-3,17){$a$} \put(0,18){\oval(24,10)[r]}
\put(30,18){\oval(24,10)[l]} \put(15,22){\oval(20,10)[b]}
\put(15,24){\oval(20,10)[t]} \put(12,17){\circle*{1.5}}
\put(18,17){\circle*{1.5}}

\put(40,5){\circle*{1.7}} \put(40,12){\vector(1,0){3}}
\put(65,12){\vector(1,0){3}} \put(56.5,30){\vector(-1,0){3}}
\put(40,5){\vector(1,0){30}} \put(70,5){\vector(0,1){30}} \put(46,0){$b$}
\put(71,17){$a$} \put(70,35){\vector(-1,0){30}} \put(40,5){\vector(0,1){30}}
\put(62,37){$b$} \put(37,17){$a$} \put(40,19){\oval(24,10)[tr]}
\put(40,17){\oval(24,10)[br]} \put(70,17){\oval(24,10)[bl]}
\put(70,19){\oval(24,10)[tl]} \put(55,23){\oval(20,10)[b]}
\put(55,25){\oval(20,10)[t]}

\put(80,5){\circle*{1.7}} \put(80,12){\vector(1,0){3}}
\put(105,12){\vector(1,0){3}} \put(96.5,30){\vector(-1,0){3}}
\put(80,5){\vector(1,0){30}} \put(110,5){\vector(0,1){30}} \put(86,0){$b$}
\put(111,17){$a$} \put(110,35){\vector(-1,0){30}} \put(80,5){\vector(0,1){30}}
\put(102,37){$b$} \put(77,17){$a$} \put(80,19){\oval(24,10)[tr]}
\put(80,17){\oval(24,10)[br]} \put(110,18){\oval(24,12)[bl]}
\put(110,18){\oval(24,12)[tl]} \put(95,23){\oval(20,10)[b]}
\put(95,25){\oval(20,10)[t]} \put(98,18){\circle*{1.5}}

\end{picture}
\caption{Linking number of degree two} \label{fig888}
\end{figure}

\subsection{Example}

We describe here a second-degree invariant of links in $K^2 \times \R^1$
related with the self-intersection set of the discriminant (i.e. with the set
of links $f$ having two intersections of their two parts): more precisely, with
the component of this self-intersection set containing the map shown in Fig.
\ref{fig888} left. In particular, we will show that this invariant separates
the link with contractible parts shown in the center of Fig.~\ref{fig888} from
the trivial one.
\medskip

The stratum $\Delta$ of the resolved discriminant $\sigma$, related with this
component of the self-intersection set, lies in the term $\sigma_2 \setminus
\sigma_1$ of our filtration and consists of the points encoded by the data
$(\langle (x,y),(x',y')\rangle , f , \lambda)$, where
\begin{itemize}
\item
$x$ and $x'$ are some points of $S^1_1$, $y$ and $y'$ some points of $S^1_2$
(the pairs $(x,y)$ and $(x',y')$ should be different points of $S^1_1 \times
S^1_2$, although $x$ can be equal to $x'$ or $y$ can be equal to $y'$);

\item
$f \in C^\infty({\mathcal C},M^3)$, and $f(x)=f(y)$, $f(x') = f(y')$. The knots
$f(S^1_1)$ and $f(S^1_2)$ are contractible in $K^2 \times \R^1$; for any path
$s$ from $x$ to $x'$ in $S^1_1$ and any path $s'$ from $y$ to $y'$ in $S^1_2$,
the loop $f(s) \cup f(s')$ is freely homotopic to $b$ (and to $b^{-1}$), see
Fig. \ref{fig888} left;

\item
$\lambda$ is some point of the open interval, whose missing endpoints are
related with the pairs $( x,y )$ and $( x',y')$. This interval consists of two
coherent segments from the construction of the simplicial resolution,
connecting the point of $Y(2)$ corresponding to the collection $\langle (x,y),
(x',y')\rangle$ (if all points $x, y, x', y'$ are different) or $\langle x=x',
y, y' \rangle$ or $\langle x, x', y=y' \rangle$ (if some two points coincide)
with two points of $Y(1)$ corresponding to the subordinate collections $\langle
x,y\rangle$ and $\langle x',y' \rangle$. The latter endpoints of these segments
belong to the lower term $\sigma_1$ of our filtration and therefore do not
belong to our stratum of $\sigma_2 \setminus \sigma_1$.
\end{itemize}

This stratum $\Delta$ is obviously homeomorphic to a smooth manifold. (More
precisely, such strata with $f \in {\mathcal F}_d$ are homeomorphic to smooth
manifolds if $d$ is large enough). Moreover, it is orientable. Indeed, if we
walk in it in such a way that both points
\begin{equation}
f(x) \equiv f(y) \qquad \mbox{and} \qquad f(x') \equiv f(y') \label{cru}
\end{equation}
stay in a simply-connected domain in $K^2 \times \R^1$ (maybe meeting and/or
permuting there), then its canonical orientation is derived from a chosen
orientation of this domain in the same way as for links in any orientable
manifold. If we go along a large loop in our stratum, then the orientation of
the stratum is preserved along this loop if and only if the continuations of
the chosen orientation of this simply-connected domain along the traces of the
points (\ref{cru}) in $K^2 \times \R^1$ are either both preserved or both
reversed (i.e. the union of these traces defines an orientation-preserving
element of $H_1(K^2, \Z_2)$). This condition is satisfied for any loop starting
and finishing at one and the same point of our stratum. Indeed, such a loop
brings the link $f$ to its initial position, therefore the numbers of rotations
of these two points (\ref{cru}) along the disorienting direction $\{a\}$ should
be equal to one another, and the sum of these numbers is even.

So, the fundamental cycle of this stratum $\Delta$ defines a non-zero element
in the top-dimensional integer homology group of locally finite chains in
$\sigma_2 \setminus \sigma_1$.

Further, we need to prove that this cycle can be extended to a cycle in entire
$\sigma_2$, i.e. its boundary in $\sigma_1$ is zero-homologous there. To prove
it, it is enough to find a set of generators of the group $H_1(\sigma_1)$ (or
$\pi_1(\sigma_1)$) whose intersection indices with this boundary $\partial
\Delta$ are equal to 0. Let us do it.

The term $\sigma_1$ is the space of pairs $(\langle x,y \rangle, f)$, where $x
\in S^1_1$, $y \in S^1_2$, $f \in C^\infty({\mathcal C},M^3)$, $f(x)=f(y)$, and
both loops $f(S^1_1), f(S^1_2)$ are contractible. The boundary in it of the
stratum $\Delta$ consists of some such pairs with $f$ having an additional
intersection of $f(S^1_1)$ and $f(S^1_2)$. The homotopy exact sequence of the
fiber bundle $p:\sigma_1 \to M^3$, $p(\langle x,y \rangle , f) =f(x),$ contains
the fragment
\begin{equation} \label{exact}
\dots \to \Z^2 \to \pi_1(\sigma_1) \to \pi_1(M^3) \to \dots \ : \end{equation}
indeed, the fundamental group of the fiber of this bundle is generated by the
cyclic reparameterizations of $S^1_1$ and $S^1_2$ and has no additional
elements because $\pi_2(M^3)=0$. Obviously, such generators do not meet
$\partial \Delta$. Further, any element of $\pi_1(M^3)$ can be lifted to a loop
in $\sigma_1$ any point of which is a pair $(\langle x,y \rangle, f)$ such that
$f({\mathcal C})$ is contained in some small ball. These loops also do not
intersect the chain $\partial \Delta$. Therefore $\partial \Delta$ is indeed
zero-homologous in $\sigma_1$, i.e. there is a chain $\nabla$ in $\sigma_1$
such that $\partial \nabla =- \partial \Delta$. Then $\Delta+\nabla$ is a cycle
in $\sigma_2 \subset \sigma$, and its projection to $\Sigma$ is a cycle there.
By Proposition \ref{prop55} the cohomology class in $C^\infty({\mathcal
C},M^3)$ defined by the intersection numbers with this cycle is equal to zero.
Therefore the linking number with this cycle is well-defined as a homotopy
invariant of links in $C^\infty({\mathcal C},M^3)$. Let us prove that this
invariant takes non-zero value on the link shown in the center of
Fig.~\ref{fig888}.

A path in $C^\infty({\mathcal C},M^3)$ connecting this link with the trivial
one intersects the discriminant twice: first at the point shown in Fig.
\ref{fig888} right, and then at a \begin{picture}(9,7.5) \put(3,2){\oval(6,6)}
\put(6,4){\oval(6,4)[b]} \put(6,6){\oval(6,4)[t]} \put(6,2){\circle*{1.5}}
\put(9,4){\line(0,1){2}}
\end{picture}-like curve contained in some simply-connected domain. The value
of our invariant at the initial link is equal to the sum of coefficients, with
which the neighbourhoods of these points in $\Sigma$ participate in our chain
$\nabla$; these coefficients should be taken with signs $\pm$ depending on the
direction in which our path crosses $\Sigma$ with respect to its coorientation.
For the trivial discriminant point
\begin{picture}(9,7.5) \put(3,2){\oval(6,6)} \put(6,4){\oval(6,4)[b]}
\put(6,6){\oval(6,4)[t]} \put(6,2){\circle*{1.5}} \put(9,4){\line(0,1){2}}
\end{picture}
this coefficient is equal to 0, since this point can be connected with itself
in the smooth part of $\Sigma$ by a loop reversing the orientation of this
smooth part. Also, the difference between these two coefficients is equal to
$\pm 1$ (depending on the choice of the orientation of $\Delta$), because these
two intersection points are separated in $\sigma_1$ by the cycle $\partial
\Delta$. \hfill $\Box$

\end{document}